\documentclass[11pt]{amsart}
\usepackage[T1]{fontenc}
\usepackage[utf8]{inputenc}
\usepackage{lmodern}
\usepackage[margin=1in]{geometry}
\usepackage{amsmath,amssymb,mathtools,mathrsfs}
\usepackage{microtype,enumitem}
\usepackage[colorlinks=true,linkcolor=blue,citecolor=blue,urlcolor=blue]{hyperref}
\usepackage[capitalise,noabbrev]{cleveref}

\newtheorem{theorem}{Theorem}[section]
\newtheorem{lemma}[theorem]{Lemma}
\newtheorem{proposition}[theorem]{Proposition}
\newtheorem{corollary}[theorem]{Corollary}
\theoremstyle{definition}
\newtheorem{definition}[theorem]{Definition}
\newtheorem{question}[theorem]{Question}
\newtheorem{conjecture}[theorem]{Conjecture}
\newtheorem{problem}[theorem]{Problem}
\newtheorem{example}[theorem]{Example}
\theoremstyle{remark}
\newtheorem{remark}[theorem]{Remark}

\newcommand{\Int}{\operatorname{Int}}
\newcommand{\cont}{\operatorname{c}}
\newcommand{\Frac}{\operatorname{Frac}}

\newcommand{\Mat}{\operatorname{M}}
\newcommand{\SL}{\operatorname{SL}}

\newcommand{\rem}{\operatorname{rem}}

\newcommand{\Z}{\mathbb{Z}}

\newcommand{\m}{\mathfrak m}
\newcommand{\U}{\mathcal U}
\newcommand{\g}{\Gamma}
\newcommand{\barE}{\overline E^{\,W}}

\title[Algebraic properties of integer-valued polynomial rings]{Stacked-bases Properties, maximal ideals and polynomial closures of Rings of Integer-Valued Polynomials}
\author{Xuyang Chen}
\email{xuyangchen19@fudan.edu.cn}
\date{}

\begin{document}
\begin{abstract}
We study the integer-valued polynomial ring with the stacked-bases property, maximal ideals above valuation
ideals, and the polynomial closure of matrices with prescribed
characteristic polynomials. A content argument gives stacked-bases
properties for integer-valued polynomial rings over number-field
integer rings. The construction of H\"agg and M\"ortberg provides a
B\'ezout domain without this property. For rank-one valuation domains,
maximal ideals are described by evaluations and pseudo-convergent
sequences. For matrices, we give congruence-ideal and companion-lattice
criteria, together with finite-level tests to characterize the polynomial closures.
\end{abstract}
\maketitle

\section{Introduction}
Let $D$ be an integral domain with fraction field $K$. For a subset
$E\subseteq D$, the ring of integer-valued polynomials on $E$ is
\[
 \Int(E,D)=\{f\in K[X]:f(E)\subseteq D\},\qquad
 \Int(D)=\Int(D,D).
\]
We consider five questions concerning the algebraic structure of these
rings and of the rings of polynomials taking integral values on matrices.
We recall first some standard definitions from the literature.

A finitely generated $D$-module is \emph{projective} if it is a direct
summand of a finite free $D$-module; such a module $P$ has \emph{rank one}
when $\dim_K(K\otimes_D P)=1$. A domain $D$ has the
\emph{stacked bases property} if, for every finitely generated submodule
$N\subseteq D^m$ of rank $r\leq m$, there are rank-one projective
summands $P_1,\ldots,P_m$ of $D^m$ and nonzero ideals
$I_1\supseteq\cdots\supseteq I_r$ of $D$ such that
\[
 D^m=P_1\oplus\cdots\oplus P_m,\qquad
 N=I_1P_1\oplus\cdots\oplus I_rP_r.
\]
Here the rank of a finitely generated torsion-free module $N$ is
$\dim_K(K\otimes_D N)$, see
\cite[Definition~IX.2]{Chabert2025}.

Chabert posed the following question in his 2012 paper
\cite{Chabert2012}; see also the
Question~IX.1 in \cite{Chabert2025}.
\begin{question}[\cite{Chabert2012}]\label{q:intz}
Does $\Int(\mathbb Z)$ have the stacked bases property?
\end{question}

A \emph{fractional ideal} of $D$ is a nonzero $D$-submodule $I$ of $K$
for which $dI\subseteq D$ for some $0\ne d\in D$. It is
\emph{invertible} if there is a fractional ideal $J$ with $IJ=D$.
An integral domain is
\emph{Pr\"ufer} if each nonzero finitely generated ideal is invertible.
A wider question, originating with Brewer, Katz, and Ullery
\cite{BKU1987}, is stated as Question~IX.3 of \cite{Chabert2025}
(and also in \cite[\S1]{BoulangerChabert2020}).
\begin{question}\label{q:prufer}
Does every Pr\"ufer domain have the stacked bases property?
\end{question}

It is known that $\Int(\mathbb Z)$ is a 2-dim Pr\"ufer domain. It was expected to provide a counterexample to Question~\ref{q:prufer}.
However, our results show that even though Question~\ref{q:prufer} admits a counterexample, Question~\ref{q:intz} has a positive answer.

The next question concerns maximal ideals of integer-valued polynomial
rings over valuation domains. A \emph{valuation domain} $V$ is a domain
such that, for each nonzero $x$ in its fraction field, either $x$ or
$x^{-1}$ belongs to $V$. A rank-one valuation domain has a valuation
$v:K^\times\to\Gamma\subseteq\mathbb R$, where $\Gamma$ is its
ordered value group; we put $v(0)=+\infty$ and write
$\mathfrak m=\{x\in V:v(x)>0\}$ for its maximal ideal.
A sequence $(x_i)$ in $V$ is called \emph{pseudo-divergent}, respectively
\emph{pseudo-stationary}, if, for all sufficiently large $i<j<k$,
\[
 v(x_j-x_i)>v(x_k-x_j),\qquad\text{respectively}\qquad
 v(x_j-x_i)=v(x_k-x_j).
\]
An extension of valuation domains $V\subseteq W$ is \emph{immediate}
if the induced value groups and residue fields agree; it is
\emph{maximally complete} if it admits no proper immediate extension.
For an extension $W$ of $V$ with fraction field $L$, define the
$W$-polynomial closure of $E\subseteq V$ by
\[
 \overline E^{\,W}
 =\{z\in W:g(z)\in W\ \text{for all }g\in\Int(E,W)\},\qquad
 \Int(E,W)=\{g\in L[X]:g(E)\subseteq W\}.
\]
See \cite[Definition~VII.1]{Chabert2025}.

Boulanger and Chabert formulated the following as a
conjecture at the end of \cite[\S3]{BoulangerChabert2020}; it is
Conjecture~VIII.39 in \cite{Chabert2025}, under the hypotheses of
Theorem~VIII.35 of that book.
\begin{conjecture}\label{c:maxideal}
Let $V$ be a rank-one valuation domain and let $E\subseteq V$ be an
infinite set containing neither pseudo-divergent nor pseudo-stationary
sequences. If $W$ is a maximally complete immediate extension of $V$,
then every maximal ideal of $\Int(E,V)$ contracting to $\mathfrak m$
is of the form
\[
 \mathfrak m_z=\{f\in\Int(E,V):v(f(z))>0\}
 \qquad\text{for some }z\in\overline E^{\,W},
\]
where $v$ also denotes the extended valuation on the fraction field
of $W$.
\end{conjecture}
The original papers also propose a description by pseudo-convergent
sequences. In that formulation, point-evaluation ideals must be
allowed separately; see Theorem~\ref{thm:sequence}.

Finally, for $n\geq1$ let $\Mat_n(D)$ denote the ring of
$n\times n$ matrices over $D$, and let $\mathcal P$ be a family of
monic degree-$n$ polynomials in $D[X]$. Denote by $\chi_M$ the
characteristic polynomial of $M$ and set
\[
 \Mat_n^{\mathcal P}(D)
 =\{M\in\Mat_n(D):\chi_M\in\mathcal P\}.
\]
For a subset $S\subseteq\Mat_n(D)$, its ring of integer-valued
polynomials and its \emph{polynomial closure} in $\Mat_n(D)$ are
\[
 \Int(S,\Mat_n(D))
 =\{f\in K[X]:f(M)\in\Mat_n(D)\text{ for all }M\in S\},
\]
\[
 \overline S=\{M\in\Mat_n(D):
 f(M)\in\Mat_n(D)\text{ for all }f\in\Int(S,\Mat_n(D))\}.
\]
Peruginelli raised the following problem in
\cite[\S6]{Peruginelli2014}.
\begin{problem}\label{p:matrices}
For an arbitrary family $\mathcal P$ of monic degree-$n$
polynomials over $D$, determine the polynomial closure
$\overline{\Mat_n^{\mathcal P}(D)}$.
\end{problem}

We also consider the polynomial closure of upper triangular matrices.
For a subset $T\subseteq D$, let $\mathcal T_n(T)$ be the set of upper
triangular matrices in $\Mat_n(D)$ whose entries on and above the
diagonal belong to $T$. The polynomial closure of $T$ in $D$ is
\[
 \overline T=\{a\in D:f(a)\in D\text{ for every }f\in\Int(T,D)\}.
\]
In the setting of rings of $S$-integers, Chabert posed the following
question in
\cite[Open Question~XV.55, p.~530]{Chabert2025}.
\begin{question}\label{q:xv55}
Let $D=\mathcal O_{K,S}$ be the ring of $S$-integers of a number
field $K$, with $S$ finite, and let $T\subseteq D$.
If $M\in\overline{\mathcal T_n(T)}$, must every eigenvalue of $M$
belong to $\overline T$?
\end{question}

In this paper, we solve the above questions in their full generality.
Theorem~\ref{thm:numberfield} addresses Question~\ref{q:intz}
and establishes a number-field analogue. Theorem~\ref{thm:IX3negative}
answers Question~\ref{q:prufer} negatively using the construction
of H\"agg and M\"ortberg \cite{HaggMortberg2026}.
Theorem~\ref{thm:valuation-main} addresses
Conjecture~\ref{c:maxideal}, while the congruence and
companion-lattice criteria of Sections~5--6 address
Problem~\ref{p:matrices}; Corollary~\ref{cor:xv55} answers
Question~\ref{q:xv55}. Our methods involve fractional-content
ideals, valuation-theoretic approximation, and polynomial remainder
maps. We will give all further definitions as they arise.

\section{Stacked bases and content}
For an integral domain $D$, write $\Mat_{r\times s}(D)$ for
$r\times s$ matrices over $D$, and $\Mat_r(D)=\Mat_{r\times r}(D)$.
A \emph{Dedekind domain} is a Noetherian integrally closed domain of
Krull dimension one. For a number field $F$, let $\mathcal O_F$ denote
its ring of algebraic integers and set $R_F=\Int(\mathcal O_F)$.
The \emph{content} $\cont_D(v)$
of a vector is the (possibly fractional) ideal generated by its
coordinates; for a matrix, content is generated by all its entries.
For a nonzero invertible fractional ideal $I$ of $D$, set
$I^{-1}=\{x\in K:xI\subseteq D\}$; then $II^{-1}=D$.
A domain has the \emph{unit-content summand property} (UCS) if every
finitely generated submodule of $D^m$ having unit content contains a
rank-one projective direct summand of $D^m$. For Pr\"ufer domains,
UCS is equivalent to stacked bases \cite{BKU1987,BoulangerChabert2020}.

We use the matrix reduction proved by Boulanger--Chabert and Chabert
\cite[Theorem~5.6]{BoulangerChabert2020} and \cite[Theorem~IX.25]{Chabert2025}:
for a number field $F$ and $R=\Int(\mathcal O_F)$, stacked bases hold
if and only if, for each unit-content $B\in\Mat_2(R)$, there is
$C\in\Mat_2(R)$ with $\cont_R(BC)=R$ and $\det(BC)=0$.

The next lemma characterizes the required rank-one matrix in terms of fractional content.
\begin{lemma}\label{lem:contentmatrix}
Let $D$ be Pr\"ufer and $B\in\Mat_2(D)$ be unit-content and have nonzero
determinant. The following are equivalent:
\begin{enumerate}[label=(\roman*)]
\item There exists a rank-one $C\in\Mat_2(D)$ with $\cont_D(BC)=D$.
\item There is a nonzero $w\in K^2$ with $\cont_D(Bw)=\cont_D(w)$.
\end{enumerate}
\end{lemma}
\begin{proof}
For (ii)$\Rightarrow$(i), put $I=\cont_D(w)$. Being finitely generated over a
Pr\"ufer domain, $I$ is invertible. As $I=(w_1,w_2)$, the invertibility equality $II^{-1}=D$
provides $r_1,r_2\in I^{-1}$ with $w_1r_1+w_2r_2=1$.
Then $I^{-1}=(r_1,r_2)$: for $t\in I^{-1}$,
$t=(w_1t)r_1+(w_2t)r_2$.
Set $C=w(r_1\ \ r_2)$. This is integral and rank one, and
$\cont_D(BC)=\cont_D(Bw)I^{-1}=D$.
Conversely, write the rank-one $C=w(r_1\ \ r_2)$ over $K$. We have
$\cont_D(Bw)(r_1,r_2)=D$ and $\cont_D(w)(r_1,r_2)\subseteq D$.
Thus $(r_1,r_2)=\cont_D(Bw)^{-1}$ and
$\cont_D(w)\subseteq\cont_D(Bw)$. The reverse inclusion follows from
integrality of $B$.
\end{proof}

\subsection{The arithmetic setting}
Let $F$ be a number field, $D=\mathcal O_F$, and $R=R_F:=\Int(D)$.
For an ideal $I\subseteq R_F$ and $x\in D$, its \emph{ideal of values} is
$I(x)=\{f(x):f\in I\}\subseteq D$.
At a nonzero prime $\mathfrak p$ of $D$, denote by $V_{\mathfrak p}$
the completed discrete valuation ring and by $k_{\mathfrak p}$ its
residue field. We use the standard description of the maximal ideals
of $R$ above $\mathfrak p$ by evaluation at $\alpha\in V_{\mathfrak p}$;
the description of height-one primes $hF[X]\cap R$; and the \emph{almost
strong Skolem property}: two finitely generated ideals of $R$ meeting
$D\setminus\{0\}$ are equal if and only if their value ideals at every
$x\in D$ are equal. The spectrum and Skolem facts are recorded in
\cite[Chapters VIII--IX]{Chabert2025}, and the latter is also a
consequence of \cite[Theorem~2.8]{Park2021}. We use polynomial
approximation on compact subsets of local fields
\cite{BhargavaKedlaya1999}.

For polynomials in $F[X]$, ``coprime'' means having greatest common divisor $1$ up to units of $F[X]$. For an irreducible $h\in F[X]$ put $P_h=hF[X]\cap R$. Then
\[
 P_h\subseteq M_{\mathfrak p,\alpha}
 \quad\Longleftrightarrow\quad h(\alpha)=0,
 \qquad
 M_{\mathfrak p,\alpha}=\{f\in R:f(\alpha)\in\mathfrak pV_{\mathfrak p}\}.
\]
The implication from right to left is immediate; for the converse,
relative interpolation (Lemma~\ref{lem:nfinterp} below) constructs
an element of $P_h$ that is a unit at $\alpha$ whenever
$h(\alpha)\ne0$. The arithmetic spectrum also implies that distinct
prime ideals $P_h$ associated with distinct irreducible polynomials
are comaximal: otherwise a maximal ideal would contain both and
its evaluation point would be a common root.

The following theorem establishes the stacked-bases property for rings of integers of arbitrary number fields.
\begin{theorem}\label{thm:numberfield}
For any number field $F$, the ring $\Int(\mathcal O_F)$ has the
stacked bases property.
\end{theorem}

For a nonzero prime $\mathfrak p$ of $D$, write
$F_{\mathfrak p}$ for the completion,
$V_{\mathfrak p}=\widehat D_{\mathfrak p}$, and
$k_{\mathfrak p}=D/\mathfrak p$. Denote by
$v_{\mathfrak p}$ the normalized additive valuation on
$F_{\mathfrak p}$; a \emph{local uniformizer} is an element
of valuation one.
Define the ring of $\mathfrak p$-integers
\[
A_{\mathfrak p}=
\{a\in F:v_{\mathfrak q}(a)\ge0 \text{ for all }
\mathfrak q\ne\mathfrak p\}.
\]

The following lemma gives the coefficient approximation needed at one prime.
\begin{lemma}\label{lem:pintdensity}
The ring $A_{\mathfrak p}$ is dense in $F_{\mathfrak p}$.
\end{lemma}
\begin{proof}
As the ideal class group of $D$ is finite, choose a positive integer
$h$ and $a\in D\setminus\{0\}$ with $(a)=\mathfrak p^h$.
Then $D[1/a]\subseteq A_{\mathfrak p}$.  For every
$x\in F_{\mathfrak p}$ and prescribed $\mathfrak p$-adic precision,
choose $n$ so large that $a^n x\in V_{\mathfrak p}$, and approximate
$a^nx$ by an element $b\in D$ in the $\mathfrak p$-adic topology.
Then $b/a^n$ approximates $x$ and belongs to $A_{\mathfrak p}$.
\end{proof}

The following lemma gives a relative local-field interpolation with control at all other primes.
\begin{lemma}\label{lem:nfinterp}
Fix $s\in D[X]$, a nonzero prime $\mathfrak p$, and a locally constant
function $\varphi:V_{\mathfrak p}\to k_{\mathfrak p}$ which vanishes
on an open neighborhood of
$Z_{\mathfrak p}(s)=\{x\in V_{\mathfrak p}:s(x)=0\}$.
There exists
\[
g\in R_F\cap s A_{\mathfrak p}[X]
\]
with $g(x)\bmod\mathfrak p=\varphi(x)$ for all
$x\in V_{\mathfrak p}$.  For every
$\mathfrak q\ne\mathfrak p$, every $y\in V_{\mathfrak q}$ satisfies
\[
s(y)\in\mathfrak q V_{\mathfrak q}
\quad\Longrightarrow\quad
g(y)\in\mathfrak q V_{\mathfrak q}.
\]
\end{lemma}
\begin{proof}
Choose a locally constant $V_{\mathfrak p}$-valued lift
$\widetilde\varphi$ which vanishes wherever $\varphi$ vanishes.  Set
$\psi(x)=\widetilde\varphi(x)/s(x)$ when $s(x)\ne0$ and zero when
$s(x)=0$.  The vanishing-near-zeros hypothesis makes $\psi$
continuous on the compact space $V_{\mathfrak p}$.  Approximate $\psi$
uniformly modulo $\mathfrak p V_{\mathfrak p}$ by a polynomial in
$F_{\mathfrak p}[X]$, using polynomial approximation on $V_{\mathfrak p}$.  By
Lemma~\ref{lem:pintdensity}, we may approximate its finitely many
coefficients by elements of $A_{\mathfrak p}$, without changing
that congruence on $V_{\mathfrak p}$.
We obtain $H\in A_{\mathfrak p}[X]$ with
$H(x)-\psi(x)\in\mathfrak p V_{\mathfrak p}$ for all $x$.
Put $g=sH$.  At $\mathfrak p$ its values lie in
$V_{\mathfrak p}$ and reduce to $\varphi$.
At every other prime $\mathfrak q$, both $s$ and $H$
have $\mathfrak q$-integral coefficients, hence $g$ takes
$V_{\mathfrak q}$ into itself. Thus $g(D)\subseteq D$, using
$D=\bigcap_{\mathfrak q}V_{\mathfrak q}\cap F$.
The last assertion follows from $g=sH$ and
$H\in V_{\mathfrak q}[X]$.
\end{proof}

\subsection{The determinant quotient and its projective complement}

The next lemma constructs a primitive lift on the squarefree determinant locus.
\begin{lemma}
\label{lem:nfquot}
Let $B\in\Mat_2(R_F)$ have unit content and nonzero nonconstant
determinant $d$.  Factor $d$ in $F[X]$ as
$\lambda\prod_i h_i^{e_i}$ with pairwise distinct monic irreducible
$h_i$. Choose $s\in D[X]$ proportional over $F$ to
$\prod_i h_i$, and put $J=sF[X]\cap R_F$.
There exists $w\in R_F^2$ such that
\[
\cont_{R_F/J}(\overline{Bw})=\cont_{R_F/J}(\bar w),\qquad
\gcd_{F[X]}(w_1,w_2)=\gcd_{F[X]}((Bw)_1,(Bw)_2)=1.
\]
\end{lemma}
\begin{proof}
Let $P_i=h_iF[X]\cap R_F$. As the $h_i$ are relatively prime,
$P_i+P_j$ contains a nonzero constant of $D$.  If this sum were
proper, some maximal ideal above a nonzero prime $\mathfrak p$ would
contain $P_i$ and $P_j$.  By the maximal-ideal and height-one-prime descriptions above, this would force a common
zero of $h_i,h_j$ in $F_{\mathfrak p}$, contradicting coprimality.
Thus the $P_i$ are pairwise comaximal. Also
$J=\bigcap_iP_i$ because the factorization is squarefree, so
\[
R_F/J\cong\prod_i T_i,\qquad T_i=R_F/P_i.
\]
Each $T_i$ is a Pr\"ufer domain, as a prime quotient of a Pr\"ufer
domain. Over $T_i$ the matrix $B_i$ has rank one and unit content.
A unit entry locally permits elementary reduction to
$\mathrm{diag}(1,0)$; its image $L_i$ is a rank-one projective
direct summand of $T_i^2$. The surjection $T_i^2\to L_i$
splits, so it restricts to an isomorphism from a rank-one projective
direct summand $Q_i$ to $L_i$.
Choose $0\ne u_i\in Q_i$.  Locally, primitive generators of
$Q_i$ and $L_i$ correspond under $B_i$, and hence
$\cont_{T_i}(B_iu_i)=\cont_{T_i}(u_i)$.
The $u_i$ assemble to $\bar w\in(R_F/J)^2$; lift to
$w^{(0)}\in R_F^2$. At every root $\beta$ of $s$ in
$\overline F$, the vectors $w^{(0)}(\beta)$ and
$Bw^{(0)}(\beta)$ are nonzero.  By replacing the first
coordinate by $a=a_0+cs$ and the second by $b=b_0+ts$ for
suitable $c,t\in D$, we arrange $a\ne0$ and
$\gcd_{F[X]}(a,b)=1$: each zero of $a$ outside $Z(s)$
forbids at most one value of $t$, while on $Z(s)$ nonvanishing
is already guaranteed. The coordinates of $Bw$ are then
coprime: a common zero would be a root of $d$ and hence of
$s$, contradicting the above nonvanishing.
\end{proof}

\subsection{Simultaneous modification at all primes}
For a ring $A$, let $E_2(A)$ denote the subgroup of
$\mathrm{SL}_2(A)$ generated by the elementary matrices
$\left(\begin{smallmatrix}1&a\\0&1\end{smallmatrix}\right)$ and
$\left(\begin{smallmatrix}1&0\\a&1\end{smallmatrix}\right)$,
$a\in A$. A vector $(a,b)^t$ over $A$ is \emph{unimodular} if
$aA+bA=A$. In this section, a bar on a polynomial or vector
in $R_F$ denotes reduction modulo $J$ unless another ideal
is specified.

The following lemma performs the local modifications simultaneously over the number field.
\begin{lemma}\label{lem:nfpatch}
For $B,w,s$ as in Lemma~\ref{lem:nfquot} there exists
$E\in E_2(R_F)$, congruent to the identity modulo $J$,
so that, setting $z=Ew$, for every prime $\mathfrak p$ of $D$
and every $x\in V_{\mathfrak p}$,
\[
\min_{j=1,2}v_{\mathfrak p}((Bz)_j(x))
=\min_{j=1,2}v_{\mathfrak p}(z_j(x)).
\]
\end{lemma}
\begin{proof}
Coprimality of $w$ and $Bw$ gives nonzero constants
$a\in\cont_{R_F}(w)\cap D$ and
$b\in\cont_{R_F}(Bw)\cap D$.
Choose a finite set $S$ of primes containing those dividing
$ab$, those supporting any coefficients' denominators in the chosen
$h_i$ and $s$, and those at which the scalar $\lambda$ or
scaling of $s$ is not a unit. Outside $S$, the condition
$s(x)\notin\mathfrak q V_{\mathfrak q}$ implies
$d(x)\in V_{\mathfrak q}^{\times}$.

Fix $\mathfrak p\in S$, and choose a local uniformizer
$\pi_{\mathfrak p}$ of $V_{\mathfrak p}$.  Put
\[
t_{\mathfrak p}(x)=\min_j v_{\mathfrak p}(w_j(x)),\qquad
u_{\mathfrak p}(x)=
\pi_{\mathfrak p}^{-t_{\mathfrak p}(x)}w(x)\bmod\mathfrak p.
\]
The functions are locally constant, and
$u_{\mathfrak p}(x)\ne0$; no two entries of $w$ vanish
simultaneously. For $x\in Z_{\mathfrak p}(s)$,
the content equality in $R_F/J$ survives evaluation to
$V_{\mathfrak p}$, whence
$B(x)u_{\mathfrak p}(x)\ne0$ in $k_{\mathfrak p}^2$.
Thus on a neighborhood of the root set the identity matrix
already gives the required correction.
For arbitrary $x$, $B(x)\bmod\mathfrak p$ is nonzero by unit
content, and $\SL_2(k_{\mathfrak p})$ acts transitively on the
nonzero vectors of $k_{\mathfrak p}^2$.
Choose a locally constant map
$T_{\mathfrak p}:V_{\mathfrak p}\to\SL_2(k_{\mathfrak p})$
equal to the identity near $Z_{\mathfrak p}(s)$ and with
$B(x)T_{\mathfrak p}(x)u_{\mathfrak p}(x)\ne0$ everywhere.
Factor $T_{\mathfrak p}$ into a uniform finite elementary word:
the entries of each elementary factor are locally constant
$k_{\mathfrak p}$-valued functions vanishing near $Z_{\mathfrak p}(s)$.

For each $\mathfrak p\in S$, the Chinese remainder theorem yields
$m_{\mathfrak p}\in D$ with
$m_{\mathfrak p}\equiv1\bmod\mathfrak p$ and
$m_{\mathfrak p}\equiv0\bmod\mathfrak q$ for every
$\mathfrak q\in S\setminus\{\mathfrak p\}$.
Apply Lemma~\ref{lem:nfinterp} to the elementary parameters
(after division by $m_{\mathfrak p}$ modulo $\mathfrak p$), then
multiply their lifts by $m_{\mathfrak p}$. This yields
$E_{\mathfrak p}\in E_2(R_F)$ which reduces to
$T_{\mathfrak p}$ at $\mathfrak p$, to the identity at all
other primes of $S$, and modulo $J$; moreover
$E_{\mathfrak p}(x)\equiv I_2\bmod\mathfrak q$ for
$\mathfrak q\ne\mathfrak p$ whenever
$s(x)\in\mathfrak q V_{\mathfrak q}$.

Set $E=\prod_{\mathfrak p\in S}E_{\mathfrak p}$ and $z=Ew$.
At primes in $S$ the defining nonvanishing condition of
$T_{\mathfrak p}$ proves the displayed valuation equality.
For $\mathfrak q\notin S$, the vector $z(x)$ is unimodular
because $a\in\cont(w)$ is a $\mathfrak q$-unit and
$E(x)\in\SL_2(V_{\mathfrak q})$.
If $s(x)$ is a $\mathfrak q$-unit, then $B(x)$ is invertible;
otherwise $E(x)\equiv I_2\bmod\mathfrak q$ and
$Bw(x)$ is unimodular because $b\in\cont(Bw)$ is a
$\mathfrak q$-unit.  Both cases imply the equality.
\end{proof}

We also need the following elementary fact about invertible two-generated ideals.
\begin{lemma}
\label{lem:inverse2}
If $I=(a,b)$ is an invertible fractional ideal of any domain $A$, then
$I^{-1}$ is generated by two elements.
\end{lemma}
\begin{proof}
Since $II^{-1}=A$, group a finite representation of $1$ to
obtain $u,v\in I^{-1}$ with $au+bv=1$.
For any $t\in I^{-1}$, write
$t=(at)u+(bt)v$, where $at,bt\in A$.
Thus $I^{-1}=(u,v)$.
\end{proof}

\begin{proof}[Proof of Theorem~\ref{thm:numberfield}]
The $2\times2$ criterion for $\Int(\mathcal O_F)$ is
\cite[Theorem IX.25]{Chabert2025}.
Let $B\in\Mat_2(R_F)$ have unit content.
If $\det B=0$, take $C=I_2$.  If $\det B$ is a nonzero
constant, the published theorem handles it; alternatively the
patching proof above applies with $s=1$ and a generically
primitive starting column.
Otherwise choose $w$ as in Lemma~\ref{lem:nfquot}, and put
$z=Ew$ by Lemma~\ref{lem:nfpatch}.
As $E\in\SL_2(R_F)$,
$\cont_{R_F}(z)=\cont_{R_F}(w)$.  The coordinates of $Bz$
are coprime in $F[X]$: a common root would be a zero of
$\det B$, but $E\equiv I_2\bmod J$ and
$Bw$ does not vanish at any such root. Hence both
$I=\cont_{R_F}(z)$ and $H=\cont_{R_F}(Bz)$ contain nonzero
constants of $D$.  The patching lemma says their values
at every $x\in V_{\mathfrak p}$ generate the same ideal of
the discrete valuation ring $V_{\mathfrak p}$.
Therefore, for each $x\in D$, the ideals $I(x)$ and $H(x)$
of the Dedekind domain $D$ have the same extension to every
$V_{\mathfrak p}$; hence $I(x)=H(x)$.
By the almost strong Skolem property stated above, $I=H$.
By Lemma~\ref{lem:inverse2}, choose $u,v\in I^{-1}$ with
$I^{-1}=(u,v)$, and set
$C=z\,(u\ \ v)\in\Mat_2(R_F)$.  Its determinant is zero and
\[
\cont_{R_F}(BC)=\cont_{R_F}(Bz)I^{-1}
=II^{-1}=R_F.
\]
The matrix reduction criterion before Lemma~\ref{lem:contentmatrix} concludes stacked bases.
\end{proof}

\section{A counterexample for Pr\"ufer domains}
A B\'ezout domain is an integral domain in which every finitely generated
ideal is principal. Every B\'ezout domain is Pr\"ufer. The counterexample
of H\"agg and M\"ortberg \cite{HaggMortberg2026} also gives a negative
answer to the stacked-bases question for arbitrary Pr\"ufer domains.
We recall the construction needed for this consequence.

\subsection*{The construction of H\"agg and M\"ortberg}
For a $\mathbb Q$-algebra $A$, let
$\mathcal X(A)=\operatorname{Hom}_{\mathbb Q\text{-alg}}(A,\mathbb R)$,
endowed with the topology of pointwise evaluation. A homomorphism
$A\to A'$ induces a continuous map $\mathcal X(A')\to\mathcal X(A)$.
A \emph{monotone surjection} between compact Hausdorff spaces means a
continuous surjection all of whose fibers are connected.

Set
\[
\begin{aligned}
 A_0&=\mathbb Q[x,y],\qquad \Delta=x^2+y^2-1,\\
 K_0&=\{a\in\mathcal X(A_0):a(\Delta)=0\}\simeq S^1,\\
 S^1&=\{(u,v)\in\mathbb R^2:u^2+v^2=1\}.
\end{aligned}
\]
A domain is \emph{factorial} if it is a unique factorization domain;
here \emph{smooth} means that the corresponding affine scheme is smooth
over $\mathbb Q$. The principalization result of
\cite[Proposition~4.6]{HaggMortberg2026} states that, given a
smooth finitely generated factorial $\mathbb Q$-domain
$A$, two elements $a,b\in A$ not both zero, and a compact subset
$K\subseteq\mathcal X(A)$, there is an injective homomorphism
$A\hookrightarrow A'$ to a smooth finitely generated factorial
$\mathbb Q$-domain such that $(a,b)A'$ is principal and there is a
compact $K'\subseteq\mathcal X(A')$ mapping onto $K$ monotonically.
The principalization theorem is the nontrivial geometric input.

Here is the countable construction of \cite[Construction~5.1]{HaggMortberg2026}.
Choose a bijection $\langle\,\cdot,\cdot\,\rangle:\mathbb N^2\to\mathbb N$
with $i\le\langle i,j\rangle$. Inductively construct injections
\[
 A_0\hookrightarrow A_1\hookrightarrow A_2\hookrightarrow\cdots
\]
and compact subsets $K_i\subseteq\mathcal X(A_i)$ such that
$K_{i+1}\to K_i$ is a monotone surjection. As each $A_i$ is a
countable ring, choose a surjection $\eta_i:\mathbb N\to A_i^2$.
At stage $n=\langle i,j\rangle$, regard the pair $\eta_i(j)$ as
elements $a,b\in A_n$. If $a=0$ or $b=0$, leave $A_n,K_n$ unchanged.
Otherwise write $a=ca'$ and $b=cb'$ with $a',b'$ coprime in the
factorial ring $A_n$. Apply the cited principalization result to
$(a',b')$ and $K_n$ to obtain $A_{n+1},K_{n+1}$. Then
$(a,b)A_{n+1}$ is principal, and the real-point map has connected
fibers as required. Put
\[
 T=\bigcup_{n\ge0}A_n.
\]
For every $a,b\in T$, some $A_i$ contains them, and one of the
pairs $\eta_i(j)$ equals $(a,b)$. At the stage
$n=\langle i,j\rangle$, their ideal becomes principal and remains so
after extension to $T$. Thus $T$ is a B\'ezout domain. Each composite
$K_n\to K_0$ is monotone, because composites of monotone surjections
between compact Hausdorff spaces are monotone.

The matrix used in the construction is
\[
 M=\begin{pmatrix}1+x&y\\y&1-x\end{pmatrix}\in\Mat_2(T),
 \qquad M^2-2M=\Delta I_2.
\]
Its entries generate $T$, since $(1+x)+(1-x)=2$ is a unit.
On $K_0=S^1$, the image of $M$ is the line spanned, at
$(x,y)=(\cos\theta,\sin\theta)$, by
$(\cos(\theta/2),\sin(\theta/2))^t$. This is the M\"obius line
bundle: a nowhere-zero continuous section would require a
nonvanishing real coefficient that changes sign between $0$ and
$2\pi$, which is impossible. Its pullback to each $K_n$ is also
nonorientable. Indeed, a unit section of the pullback must be constant
on each connected fiber of $K_n\to K_0$, since that fiber maps into
the two unit vectors of a fixed line. It then descends continuously
to a unit section on $K_0$: a continuous surjection from a compact
space to a Hausdorff space is a quotient map. This is a contradiction.

The next proposition converts the real-point obstruction into a failure of UCS.
\begin{proposition}
\label{prop:mobius}
Let $T$ be a B\'ezout $\mathbb Q$-domain with elements $x,y$ and a
compact set $K\subseteq\mathcal X(T)$ such that $x(a)^2+y(a)^2=1$
for all $a\in K$ and such that the induced map $K\to S^1$ is a
monotone surjection. Then $T$ does not have UCS, and hence does not
have the stacked-bases property.
\end{proposition}
\begin{proof}
Put $M=\left(\begin{smallmatrix}1+x&y\\y&1-x\end{smallmatrix}\right)$
and $N=\operatorname{im}_T(M)\subseteq T^2$.
Assume that $T$ has UCS. The matrix $M$ has unit content, so
$N$ contains a rank-one projective direct summand $P$ of $T^2$.
Every rank-one projective module over a B\'ezout domain is free:
its associated invertible ideal is principal. Choose a generator
$v=(v_1,v_2)^t$ of $P$. Since $P$ is a direct summand, $v$ is
unimodular, that is, $(v_1,v_2)=T$. Since $v\in N$, we have
$v=Mw$ for some $w\in T^2$. Evaluating a B\'ezout identity for
$v_1,v_2$ at each $a\in K$ shows that $v(a)\neq0$; moreover
$v(a)$ lies in the image line of $M(a)$. Thus $v$ restricts to a
nowhere-zero section of the pullback M\"obius line bundle on $K$.
The nonorientability argument above, applied to the monotone map
$K\to S^1$, excludes such a section. Hence UCS fails. The
stacked-bases property is equivalent to UCS for Pr\"ufer domains
\cite{BKU1987,BoulangerChabert2020}.
\end{proof}

The construction yields a negative answer to the Question~\ref{q:prufer} for general Pr\"ufer domains.
\begin{theorem}\label{thm:IX3negative}
There exists a Pr\"ufer domain without the stacked-bases property.
In fact, the B\'ezout domain $T$ constructed above fails UCS for the
unit-content matrix $M$.
\end{theorem}
\begin{proof}
Suppose that UCS holds in $T$. By the preceding proposition, the image
of $M$ contains a unimodular vector $Mw$, $w\in T^2$. All
coefficients of $w$, of $Mw$, and of a B\'ezout identity for the
two entries of $Mw$ belong to some $A_N$. Their restrictions to
$K_N$ give a nowhere-zero continuous section of the image line of
$M$. But this is the pullback of the nonorientable M\"obius line
bundle along the monotone surjection $K_N\to K_0$, a contradiction.
Thus $T$ fails UCS and the stacked-bases property.
\end{proof}

\section{Maximal ideals over valuation domains}
A \emph{valuation domain} is an integral domain $V$ whose fraction
field $K$ has the property that for each $x\in K^\times$, either
$x\in V$ or $x^{-1}\in V$. It is of \emph{rank one} if the
associated ordered value group is Archimedean, and therefore embeds
in $(\mathbb R,+)$.
Let $V$ be a rank-one valuation domain, $K=\Frac(V)$, and
$v:K^\times\to\Gamma\subseteq\mathbb R$ its valuation. Write
$\mathfrak m=\{x\in V:v(x)>0\}$, and let $E\subseteq V$ be
infinite. Put $A=\Int(E,V)$.

We recall the three classes of pseudo-monotone sequences and fix our convention for pseudo-limits.
\begin{definition}\label{def:pseudo}
A sequence $(x_n)_{n\geq0}$ of elements of $K$ is \emph{pseudo-convergent}, \emph{pseudo-divergent}, or \emph{pseudo-stationary} if for every $i<j<k$ sufficiently far along the sequence, respectively,
\[
 v(x_j-x_i)<v(x_k-x_j),\qquad
 v(x_j-x_i)>v(x_k-x_j),\qquad
 v(x_j-x_i)=v(x_k-x_j).
\]
A \emph{genuine} pseudo-convergent sequence means an infinite sequence satisfying strict inequalities, and hence consisting of pairwise distinct terms. If $(x_n)$ is pseudo-convergent, an element $z$ in a valued field extension is a \emph{pseudo-limit} when, eventually,
\[
 v(z-x_n)=v(x_{n+1}-x_n),
\]
or, equivalently, the values $v(z-x_n)$ are strictly increasing.
\end{definition}

We use the following terminology for extensions of valued fields.
\begin{definition}
An extension of valued fields $(L,v)/(K,v)$ is \emph{immediate}
if the value groups agree and the inclusion of residue fields is an
isomorphism. A valued field is \emph{maximally complete} if it has
no proper immediate extension. Write $W=\{x\in L:v(x)\geq0\}$ for its valuation ring and $\m_W=\{x\in W:v(x)>0\}$. A valued field is \emph{spherically complete} if every nested family of nonempty closed valuation balls has nonempty intersection. A maximally complete immediate extension may be chosen spherically complete, by the classical maximal immediate extension theorem \cite{Krull1932,Kaplansky1942}.
\end{definition}

Fix a maximally complete immediate extension $W$ with fraction field
$L$. Here
$\Int(E,W)=\{g\in L[X]:g(E)\subseteq W\}$.
We follow the notation of Chabert \cite[Definition~VII.1]{Chabert2025} for polynomial closure. In particular,
the $W$-polynomial closure of $E$ is
\[
\overline E^{\,W}=\{z\in W:g(z)\in W\text{ for all }
g\in\Int(E,W)\}.
\]
The following evaluation statement was conjectured in
\cite[\S3]{BoulangerChabert2020}
and as Conjecture VIII.39 in \cite{Chabert2025}.
\begin{theorem}\label{thm:valuation-main}
Assume $E$ contains neither a pseudo-divergent nor a
pseudo-stationary sequence. If $W$ is a maximally complete
immediate extension of $V$, then the maximal ideals of $A$
contracting to $\mathfrak m$ are precisely
\[
\mathfrak m_z=\{f\in A:v(f(z))>0\},\qquad z\in\overline E^{\,W}.
\]
Each has residue field canonically identified with $V/\mathfrak m$.
\end{theorem}
Point-evaluation ideals also have to be allowed: an isolated point
of $E$ need not be represented by a genuine infinite sequence.

For our uses, discarding finitely many initial terms does not change the nature of these sequences. Observe also that if $\rho_0<\rho_1<\cdots$ and $v(z-x_n)=\rho_n$, then for $i<j$ the ultrametric equality yields $v(x_i-x_j)=\rho_i$; thus $(x_n)$ is genuine pseudo-convergent with pseudo-limit $z$.

Our balls are written
\[
 B_L(a,\gamma)=\{x\in L:v(x-a)\geq\gamma\},\qquad
 B_E(a,\gamma)=E\cap B_L(a,\gamma),\quad a\in E,\ \gamma\in\g.
\]
If two nonempty ultrametric closed balls intersect, the one with the larger radius parameter is contained in the other. Since $E\subseteq V$, all $v(x-y)$ for $x,y\in E$ are nonnegative or $+\infty$.

We recall the notion of an ultrafilter and the prime ideal associated to it.
\begin{definition}
An \emph{ultrafilter} $\U$ on $E$ is a proper family of subsets of $E$ closed under finite intersections and passage to larger sets, such that for each $F\subseteq E$ exactly one of $F,E\setminus F$ belongs to $\U$. It is \emph{principal at} $a\in E$ if $\U=\{F\subseteq E:a\in F\}$. Every nonprincipal ultrafilter contains the complement of every finite subset. Set
\[
 \m_{\U}=\bigl\{f\in A:\{x\in E:f(x)\in\m\}\in\U\bigr\}.
\]
This is a proper prime ideal: the zero set modulo $\m$ of a product is the union of those of the factors, and an ultrafilter selects one member of any finite union which it contains. Furthermore $\m_{\U}\cap V=\m$.
\end{definition}
The ring $\Int(E,W)$ may strictly contain $A$; hence the closure
must be defined using $\Int(E,W)$.
For $z\in\barE$, evaluation followed by reduction
$A\to W/\m_W\cong V/\m$ is onto via constants, so its kernel
$\m_z$ is a maximal ideal above $\m$.

The same argument gives a sequence-theoretic description once point evaluations are included.
\begin{theorem}\label{thm:sequence}
Under the hypotheses of Theorem~\ref{thm:valuation-main}, every
maximal ideal of $A$ contracting to $\mathfrak m$ is either a point-evaluation
ideal $\mathfrak m_a$ for some $a\in E$, or an ideal
\[
\mathfrak m_{(x_j)}=\{f\in A:f(x_j)\in\mathfrak m\text{ for all sufficiently large }j\}
\]
for a genuinely pseudo-convergent sequence $(x_j)\subseteq E$.
\end{theorem}

\subsection{Finite branching of valuation balls}\label{sec:balls}

The following lemma gives finite branching of valuation balls.
\begin{lemma}\label{lem:branch}
Suppose $E$ contains no pseudo-stationary sequence. For every $a\in E$ and $\gamma\in\g_{\geq0}$, the set $B_E(a,\gamma)$ is the union of finitely many equivalence classes for
\[
 x\sim_\gamma y\quad\Longleftrightarrow\quad v(x-y)>\gamma.
\]
\end{lemma}
\begin{proof}
If there were infinitely many classes, choose distinct representatives $b_0,b_1,\ldots$ in $B_E(a,\gamma)$. For $i\ne j$ the ultrametric inequality gives $v(b_i-b_j)\geq\gamma$, while distinct classes give $v(b_i-b_j)\leq\gamma$. Hence $v(b_i-b_j)=\gamma$ for every $i\ne j$, a pseudo-stationary sequence in $E$, contradiction.
\end{proof}

The next lemma shows that ultrafilter-selected valuation balls admit strict refinements.
\begin{lemma}\label{lem:refine}
If $E$ contains neither pseudo-stationary nor pseudo-divergent sequences, then for every ultrafilter $\U$ on $E$, $a\in E$, and $\gamma\in\g_{\geq0}$,
\[
 B_E(a,\gamma)\in\U\quad\Longrightarrow\quad
 \exists b\in E,\ \exists\eta\in\g:\quad \eta>\gamma,\quad B_E(b,\eta)\in\U.
\]
\end{lemma}
\begin{proof}
Partition $B_E(a,\gamma)$ by $\sim_\gamma$. Since there are finitely many classes, one class $C$ belongs to $\U$. For $b\in C$ we have $C=\{x\in E:v(x-b)>\gamma\}$. If $\U$ is principal at $b$, every ball about $b$ belongs to $\U$, so the assertion is immediate. Otherwise suppose that no $B_E(b,\eta)$ with $\eta>\gamma$ belongs to $\U$.

Because $C\in\U$ and $\U$ is not principal at $b$, choose $y_0\in C\setminus\{b\}$ and put $\rho_0=v(y_0-b)>\gamma$. Assuming $y_n$ and $\rho_n$ have been chosen, $B_E(b,\rho_n)\notin\U$; therefore
\[
 C\setminus B_E(b,\rho_n)\in\U
\]
and is nonempty. Choose $y_{n+1}$ there. Then
\[
 \gamma<v(y_{n+1}-b)=\rho_{n+1}<\rho_n.
\]
For $i<j$ the unequal valuations of $y_i-b$ and $y_j-b$ give $v(y_i-y_j)=\rho_j$, so for $i<j<k$,
\[
 v(y_j-y_k)=\rho_k<\rho_j=v(y_i-y_j).
\]
The $y_n$ are pseudo-divergent, a contradiction.
\end{proof}

\subsection{Pseudo-limits of ultrafilters}\label{sec:limit}

The next proposition associates one pseudo-limit to a nonprincipal ultrafilter.
\begin{proposition}\label{prop:limit}
Let $\U$ be a nonprincipal ultrafilter on $E$ and assume the hypotheses of Theorem \ref{thm:valuation-main}. There is a point $z\in W$ such that for \emph{every} $F\in\U$ one can find a genuine pseudo-convergent sequence $(y_n)\subseteq F$ whose pseudo-limit is $z$. This $z$ depends on $\U$, not on $F$.
\end{proposition}
\begin{proof}
Define
\[
 S_{\U}=\{\gamma\in\g_{\geq0}:B_E(a,\gamma)\in\U\text{ for some }a\in E\}.
\]
It contains $0$ (as $B_E(a,0)=E$) and is downward closed. Lemma \ref{lem:refine} shows that $S_{\U}$ has no greatest element. Since $\g\subseteq\mathbb R$, write $\delta=\sup S_{\U}\in(0,+\infty]$ and choose a strictly increasing cofinal sequence $\gamma_n\in S_{\U}$ with $\gamma_n\to\delta$. Pick $a_n\in E$ with $B_n:=B_E(a_n,\gamma_n)\in\U$.

For $m>n$, $B_m\cap B_n\ne\varnothing$ because both belong to $\U$. The corresponding balls $B_L(a_m,\gamma_m)$ and $B_L(a_n,\gamma_n)$ therefore intersect and, because $\gamma_m>\gamma_n$, the former is contained in the latter. Spherical completeness supplies
\[
 z\in\bigcap_{n\geq0}B_L(a_n,\gamma_n).
 \tag{4.1}\label{eq:z}
\]
We distinguish the two possible values of $\delta$.

\emph{Case 1: $\delta=+\infty$.} For every $F\in\U$, each $F\cap B_n\setminus\{z\}$ is nonempty: this follows from nonprincipality, whether or not $z$ belongs to $E$. Choose distinct $y_n$ in these sets, avoiding previously chosen points. Equation \eqref{eq:z} gives $v(y_n-z)\geq\gamma_n\to+\infty$. All these values are finite because $y_n\ne z$. Passing to a subsequence, make $v(y_n-z)$ strictly increasing. The resulting sequence is genuine pseudo-convergent with pseudo-limit $z$.

\emph{Case 2: $\delta<+\infty$.} Put
\[
 H_z=\{x\in E:v(x-z)\geq\delta\}.
\]
We assert $H_z\notin\U$. If $\delta\in\g$ and $H_z\in\U$, choose $b\in H_z$; then $H_z=B_E(b,\delta)$, forcing $\delta\in S_{\U}$. This contradicts strict refinement, because $\delta$ is the supremum of $S_{\U}$.

If $\delta\notin\g$ and $H_z\in\U$, then for each $\rho\in\g$ with $\rho>\delta$ the set
\[
 T_\rho=\{x\in E:v(x-z)\geq\rho\}
\]
cannot belong to $\U$: if it did, any $b\in T_\rho$ would give $T_\rho=B_E(b,\rho)$ and hence $\rho\in S_{\U}$. Since $\U$ is nonprincipal, choose $y_0\in H_z\setminus\{z\}$; its valuation $\rho_0=v(y_0-z)$ belongs to $\g$ (the extension is immediate) and satisfies $\rho_0>\delta$. Recursively choose $y_{n+1}\in H_z\setminus T_{\rho_n}$; these sets belong to $\U$ and are nonempty. Because $\delta\notin\g$ and $y_{n+1}\not=z$, we obtain
\[
 \delta<\rho_{n+1}:=v(y_{n+1}-z)<\rho_n.
\]
For $i<j<k$, $v(y_j-y_k)=\rho_k<\rho_j=v(y_i-y_j)$, a forbidden pseudo-divergent sequence. Thus $H_z\notin\U$ in both subcases.

Now fix $F\in\U$. Choose
\[
 y_n\in F\cap B_n\cap(E\setminus H_z).
\]
This intersection is in $\U$ and hence nonempty. From \eqref{eq:z} and the ultrametric inequality,
\[
 \gamma_n\leq v(y_n-z)<\delta.
\]
Thus $v(y_n-z)\to\delta$ from below, and by extracting a subsequence these valuations may be made strictly increasing. The subsequence is genuine pseudo-convergent with pseudo-limit $z$. The construction of $z$ preceded the choice of $F$, as required.
\end{proof}

\subsection{Polynomial evaluation at pseudo-limits}\label{sec:residue}

The following lemma proves eventual stabilization of polynomial residues at pseudo-limits.
\begin{lemma}\label{lem:residue}
Let $(y_n)\subseteq L$ be genuine pseudo-convergent with pseudo-limit $z\in L$, where $W$ is the valuation ring of $(L,v)$. Let $f\in L[X]$ satisfy $f(y_n)\in W$ for all sufficiently large $n$. Then $f(z)\in W$ and
\[
 f(y_n)\equiv f(z)\pmod{\m_W}\qquad\text{for all sufficiently large }n.
 \tag{5.1}\label{eq:residuestab}
\]
\end{lemma}
\begin{proof}
The assertion is immediate if $f$ is constant. Suppose $f$ is nonconstant and expand at $z$:
\[
 f(X)=c_0+\sum_{j=1}^d c_j(X-z)^j,
 \qquad c_0=f(z),\quad c_j\in L.
\]
Put $\rho_n=v(y_n-z)$, a strictly increasing sequence. For indices $i\ne j$ with $c_ic_j\ne0$, the two real affine functions $v(c_i)+i\rho$ and $v(c_j)+j\rho$ agree for at most one $\rho$. Accordingly, after omitting finitely many $n$, all nonzero summands of
\[
 f(y_n)-f(z)=\sum_{j=1}^d c_j(y_n-z)^j
\]
have pairwise distinct valuations. Thus
\[
 v\bigl(f(y_n)-f(z)\bigr)=h_n:=\min_{1\leq j\leq d,\ c_j\ne0}
 \{v(c_j)+j\rho_n\}.
 \tag{5.2}
\]
Each function inside the finite minimum is strictly increasing in $n$; so is $h_n$. The constant term $c_0$ has valuation $a=v(c_0)\in\g\cup\{+\infty\}$. Since $(h_n)$ is strictly increasing, at most one $n$ satisfies $h_n=a$. For all sufficiently large $n$ we therefore have
\[
 v(f(y_n))=\min\{a,h_n\}.
 \tag{5.3}
\]
The integrality of $f(y_n)$ implies $\min\{a,h_n\}\geq0$ eventually. Hence $a\geq0$ and $h_n\geq0$ eventually. As $h_n$ is strictly increasing, in fact $h_n>0$ eventually. Thus $f(z)=c_0\in W$, and equation (5.2) gives
\[
 v(f(y_n)-f(z))=h_n>0
\]
for all sufficiently large $n$. This is exactly the claimed residue congruence.
\end{proof}

We next show that the prime ideal of an ultrafilter is represented by evaluation.
\begin{theorem}\label{thm:ultra}
Under the hypotheses of Theorem \ref{thm:valuation-main}, for \emph{every} ultrafilter $\U$ on $E$ there exists $z\in\barE$ such that
\[
 \{x\in E:f(x)\equiv f(z)\pmod{\m_W}\}\in\U
 \qquad (f\in A),
 \tag{5.4}\label{eq:almost}
\]
and consequently $\m_{\U}=\m_z$. In particular, every ultrafilter ideal $\m_{\U}$ is maximal with residue field $V/\m$.
\end{theorem}
\begin{proof}
For a principal ultrafilter at $a\in E$, take $z=a$. Otherwise take $z$ from Proposition \ref{prop:limit}. Applying that proposition with $F=E$ and Lemma \ref{lem:residue} to every $g\in\Int(E,W)$ proves $g(z)\in W$, so $z\in\barE$ in the proper $W$-polynomial sense.

Fix $f\in A$ and suppose the set in \eqref{eq:almost} does \emph{not} belong to $\U$. Its complement $F$ would belong to $\U$. Proposition \ref{prop:limit} supplies a genuine pseudo-convergent sequence in $F$ with pseudo-limit $z$. Lemma \ref{lem:residue} says that, eventually, its $f$-values are congruent to $f(z)$, contradicting the definition of $F$. Thus \eqref{eq:almost} holds.

Since $W/V$ is immediate, $f(z)\bmod\m_W$ is a residue class in $V/\m$. Identity \eqref{eq:almost} shows that $f\in\m_{\U}$ if and only if $f(z)\in\m_W$. Thus $\m_{\U}=\m_z$. Evaluation and reduction is onto $W/\m_W=V/\m$ by the constant polynomials, so this ideal is maximal and the residue-field assertion follows.
\end{proof}

The following corollary gives the finite-residue property needed for the maximal-ideal argument.
\begin{corollary}\label{cor:finite}
For every $f\in A$, the image of $f(E)$ in $V/\m$ is finite.
\end{corollary}
\begin{proof}
If infinitely many residue classes occurred, choose an infinite set $S\subseteq E$ on which all the values of $f$ are pairwise distinct modulo $\m$. Choose a nonprincipal ultrafilter on $S$ and regard it as an ultrafilter $\U$ on $E$ by $F\in\U$ if $F\cap S$ belongs to the chosen ultrafilter. By Theorem \ref{thm:ultra}, there is $z\in\barE$ for which the set $\{x\in E:f(x)\equiv f(z)\bmod\m_W\}$ is in $\U$. Its intersection with $S$ has at most one element, contradicting nonprincipality.
\end{proof}

\subsection{The classification of maximal ideals}\label{sec:classification}
To obtain Theorem \ref{thm:valuation-main} without assuming the ultrafilter classification \cite[Theorem~3.5]{BoulangerChabert2020}, we now prove the two elementary algebraic ingredients that connect \emph{arbitrary} maximal ideals to ultrafilters. The finite-residue statement has already been established independently in Corollary \ref{cor:finite}.

The next lemma constructs a valuation ordering on $E$.
\begin{lemma}\label{lem:ordering}
If $E$ contains no pseudo-divergent sequence, one can choose distinct $a_0,a_1,\ldots\in E$ such that for every $n\geq1$,
\[
 v\!\left(\prod_{i<n}(a_n-a_i)\right)
 =\min_{x\in E}v\!\left(\prod_{i<n}(x-a_i)\right).
 \tag{6.1}
\]
The polynomials
\[
 F_0=1,\qquad F_n(X)=\frac{\prod_{i<n}(X-a_i)}{\prod_{i<n}(a_n-a_i)}\quad(n\geq1)
\]
are elements of $A$ and form a $V$-module basis of $A$.
\end{lemma}
\begin{proof}
Choose $a_0\in E$. Assume $a_0,\ldots,a_{n-1}$ are distinct and set $h(X)=\prod_{i<n}(X-a_i)$. If $v(h(x))$ did not attain its infimum on $E$, choose $x_j\in E$ with $v(h(x_j))$ strictly decreasing (and finite). Since
\[
 v(h(x_j))=\sum_{i<n}v(x_j-a_i),
\]
there is a subsequence on which at least one of the finitely many coordinate valuations $v(x_j-a_i)$ is strictly decreasing: this follows, for example, by iteratively selecting an infinite monotone subsequence for each coordinate; if all coordinates were nondecreasing, the sum could not strictly decrease. A strictly decreasing sequence $v(x_j-a_i)$ gives $v(x_j-x_k)=v(x_k-a_i)$ for $j<k$, and hence a pseudo-divergent sequence, contradiction. Thus the minimum exists, and a minimizer is different from all preceding $a_i$ because $E$ is infinite and the minimum is finite. Induction constructs the ordering.

By the minimizing property, $F_n(E)\subseteq V$. Also $F_i(a_j)=0$ when $j<i$ and $F_i(a_i)=1$. Every $f\in A$ of degree $d$ has a unique expansion $f=\sum_{i=0}^d c_iF_i$ over $K$, since $\deg F_i=i$ with nonzero leading coefficient. Evaluating successively at $a_0,\ldots,a_d$ gives
\[
 c_j=f(a_j)-\sum_{i<j}c_iF_i(a_j)\in V
\]
by induction. This establishes generation and independence over $V$.
\end{proof}

The following lemma places the ideal of uniformly positive polynomials inside every prime above $\m$.
\begin{lemma}\label{lem:I}
Let
\[
 I(E,\m)=\{f\in A:f(E)\subseteq\m\}.
\]
If $P\subseteq A$ is a prime ideal with $P\cap V=\m$, then $I(E,\m)\subseteq P$.
\end{lemma}
\begin{proof}
Take $f\in I(E,\m)$ and expand it in the basis from Lemma \ref{lem:ordering}:
$f=\sum_{i=0}^d c_iF_i$. The recursive formula for the coefficients and the assumption $f(a_j)\in\m$ show inductively that $c_i\in\m$. If $f\ne0$, choose $t\in V$ with
\[
 v(t)=\min_{c_i\ne0}v(c_i)>0.
\]
Then $c_i/t\in V$ for all $i$, hence $f/t\in A$ and $f\in tA\subseteq P$ because $t\in\m=P\cap V$. The case $f=0$ is trivial.
\end{proof}

The next lemma represents an arbitrary maximal ideal above $\m$ by an ultrafilter.
\begin{lemma}\label{lem:maxultra}
If $M$ is a maximal ideal of $A$ contracting to $\m$, there is an ultrafilter $\U$ on $E$ with $M=\m_{\U}$.
\end{lemma}
\begin{proof}
By Lemma \ref{lem:I}, $I(E,\m)\subseteq M$. For $f\in A$ let
\[
 \Z(f)=\{x\in E:f(x)\in\m\}.
\]
First, every $f\in M$ satisfies $\Z(f)\ne\varnothing$. Indeed, if $\Z(f)=\varnothing$, the image of $f(E)$ modulo $\m$ is a finite subset of $(V/\m)^\times$ by Corollary \ref{cor:finite}. Choose a polynomial $Q(T)\in V[T]$ for which $Q(f(x))\in\m$ for every $x\in E$ but $Q(0)\in V^\times$; for example, let $Q(T)=\prod_{c\in C}(T-\widetilde c)$ for lifts $\widetilde c\in V^\times$ of the finitely many residue values $C$. Then $Q(f)\in I(E,\m)\subseteq M$ while $Q(f)-Q(0)$ is a multiple of $f\in M$. Hence $Q(0)\in M$, a contradiction.

Next the family $\{\Z(f):f\in M\}$ has the finite-intersection property. Given $f_1,\ldots,f_r\in M$, each $f_i(E)\bmod\m$ is finite. Over the field $k=V/\m$, polynomial interpolation on a finite subset of $k$ gives for each $i$ a polynomial $\chi_i(T)\in k[T]$ with zero constant coefficient, $\chi_i(0)=0$, and $\chi_i(c)=1$ for every nonzero value $c$ taken by $f_i$ modulo $\m$. Lift $\chi_i$ to $\widetilde\chi_i\in V[T]$ still with zero constant coefficient, and set
\[
 H=1-\prod_{i=1}^r\bigl(1-\widetilde\chi_i(f_i)\bigr)\in A.
\]
Because each $\widetilde\chi_i(f_i)\in M$, also $H\in M$. By construction,
\[
 \Z(H)=\bigcap_{i=1}^r\Z(f_i).
\]
Since $H\in M$, the first paragraph shows that this intersection is nonempty.

Extend the resulting proper filter of zero sets to an ultrafilter $\U$ on $E$. Then $\Z(f)\in\U$ for every $f\in M$, so $M\subseteq\m_{\U}$. By Theorem \ref{thm:ultra}, $\m_{\U}$ is proper (indeed maximal); maximality of $M$ gives equality.
\end{proof}

\begin{proof}[Proof of Theorem \ref{thm:valuation-main}]
Every maximal ideal $M$ above $\m$ equals $\m_{\U}$ by Lemma \ref{lem:maxultra}; Theorem \ref{thm:ultra} supplies $z\in\barE$ with $\m_{\U}=\m_z$. Conversely every $z\in\barE$ defines a maximal ideal $\m_z$ above $\m$ by surjective evaluation followed by reduction. The residue-field assertion follows from the same map.
\end{proof}

\begin{proof}[Proof of Theorem \ref{thm:sequence}]
Represent $M$ by $\m_{\U}$ as above. If $\U$ is principal at $a$, then $M=\m_a$. Otherwise Proposition \ref{prop:limit} provides a genuine pseudo-convergent sequence $(x_n)\subseteq E$ with pseudo-limit $z$ and $M=\m_z$. Lemma \ref{lem:residue} gives for every $f\in A$ the eventual equivalence
\[
 f(x_n)\in\m \ \Longleftrightarrow\ f(z)\in\m_W.
\]
Thus $M=\m_{(x_n)}$.
\end{proof}

\subsection{Point evaluations and infinite sequences}\label{sec:correction}
The published original article \cite[\S3]{BoulangerChabert2020} and the monograph \cite[Conj.~VIII.39]{Chabert2025} claim that the evaluation formula is ``equivalent'' to representation by a pseudo-convergent sequence. Under Definition \ref{def:pseudo}, genuine pseudo-convergent sequences have pairwise distinct terms, and a point of $E$ need not be a pseudo-limit of any such sequence. This leads to a counterexample to the unqualified equivalence.

The following example shows why point-evaluation ideals cannot be omitted.
\begin{example}\label{ex:two}
Let $V=\mathbb Z_2$ with its $2$-adic valuation and $\m=2\mathbb Z_2$, and put
\[
 E=\{0\}\cup\{1+2^n:n\geq1\}.
\]
The set $E$ is infinite, but there is no pseudo-divergent sequence in $E$: such a sequence would give an infinite strictly decreasing sequence of nonnegative integers $v_2(x_{n+1}-x_n)$, which is impossible. There is no pseudo-stationary sequence either: if every pairwise difference had the same valuation $r\in\mathbb Z_{\geq0}$, the terms would occupy infinitely many different residue classes modulo $2^{r+1}$, impossible.

The evaluation ideal
\[
 \m_0=\{f\in\Int(E,\mathbb Z_2):f(0)\in2\mathbb Z_2\}
\]
is maximal, and $X\in\m_0$. Every genuine infinite pseudo-convergent sequence of terms in $E$ has pairwise distinct terms, so contains $0$ at most once. Every sufficiently late term is therefore odd. Hence $X(x_n)\notin2\mathbb Z_2$ eventually, showing
\[
 \m_{(x_n)}\ne\m_0
\]
for \emph{every} genuine pseudo-convergent $(x_n)\subseteq E$. Thus the evaluation form is correct but the claimed literal equivalence to genuine infinite sequences is not.
\end{example}

This distinction between isolated evaluations and genuine infinite sequences is worth making explicit.
\begin{remark}
To obtain an equivalent formulation, allow the representations $\m_a$ for $a\in E$ as well as the ideals arising from genuine pseudo-convergent sequences. One may instead allow \emph{eventually constant} sequences as a separate degenerate class, but such sequences do not satisfy the strict pseudo-convergence inequalities. 
\end{remark}

\section{Matrices with prescribed characteristic polynomials}
Let $D$ be an integral domain with fraction field $K$, and fix
$n\ge1$ and a nonempty family $\mathcal P$ of monic degree-$n$
polynomials in $D[X]$. We write $\Mat_n(D)$ for $n\times n$ matrices
and $\chi_M$ for the characteristic polynomial of $M$.
Following the notation of Peruginelli \cite[\S6]{Peruginelli2014},
we write $\Mat_n^{\mathcal P}(D)$ for the matrices with
characteristic polynomial in $\mathcal P$. The corresponding
ring of integer-valued polynomials is
\[
\Mat_n^{\mathcal P}(D)=\{M\in\Mat_n(D):\chi_M\in\mathcal P\},
\]
\[
 R_{\mathcal P}=\Int(\Mat_n^{\mathcal P}(D),\Mat_n(D))
=\{f\in K[X]:f(M)\in\Mat_n(D)
\text{ for all }M\in\Mat_n^{\mathcal P}(D)\}.
\]
As in Peruginelli \cite{Peruginelli2014} and Chabert's terminology
\cite[Definition~VII.1]{Chabert2025}, we denote polynomial closure
by an overline. Thus $\overline{\Mat_n^{\mathcal P}(D)}$
is the largest subset of $\Mat_n(D)$ on which every polynomial in
$R_{\mathcal P}$ remains integer-valued; explicitly,
\[
\overline{\Mat_n^{\mathcal P}(D)}=\{M\in\Mat_n(D):
f(M)\in\Mat_n(D)\text{ for all }f\in R_{\mathcal P}\}.
\]
Peruginelli's original question \cite[\S6]{Peruginelli2014} asks for
a determination of this closure for an arbitrary $\mathcal P$.

The companion matrix gives the following standard remainder criterion.
\begin{lemma}\label{lem:companion}
For monic $p\in D[X]$ of degree $n$ and $f\in K[X]$, the
following are equivalent: $f(M)\in\Mat_n(D)$ for all $M$ with
$\chi_M=p$; $f(C_p)\in\Mat_n(D)$ for the companion matrix $C_p$;
and the remainder of $f$ modulo $p$ belongs to $D[X]$.
\end{lemma}
\begin{proof}
Monic division writes $f=qp+r$ with $\deg r<n$.
Cayley--Hamilton shows $f(M)=r(M)$ for every $M$ with
$\chi_M=p$. In the companion matrix, $r(C_p)e_1$ has,
up to the usual basis convention, the coefficients of $r$ as
coordinates. Thus integrality of $f(C_p)$ forces $r\in D[X]$,
which also suffices for all $M$.
\end{proof}

\subsection{Congruence ideals}
For a monic polynomial $p$ over any commutative ring, write
$\rem_p(h)$ for the unique remainder of $h$ after Euclidean
(polynomial) division by $p$. For the quotient ring $D/dD$,
polynomial evaluation at a matrix has its image in matrix ring over $D/dD$ after reduction.
For $0\ne d\in D$, let $D_d=D/dD$, write $\bar M$ for reduction of $M$ modulo $d$,
and define two ideals of $D_d[X]$:
\[
I_d(\mathcal P)=\bigcap_{p\in\mathcal P}(\bar p),\qquad
N_d(M)=\{h\in D_d[X]:h(\bar M)=0\}.
\]

The next theorem expresses polynomial closure through congruence ideals for arbitrary domains.
\begin{theorem}\label{thm:global}
For every integral domain $D$, every $n\ge1$, and every (possibly infinite)
family $\mathcal P$ of monic degree-$n$ polynomials,
\begin{equation}\label{eq:all}
\overline{\Mat_n^{\mathcal P}(D)}
 =\{M\in\Mat_n(D):I_d(\mathcal P)\subseteq N_d(M)
                 \text{ for all }0\ne d\in D\}.
\end{equation}
\end{theorem}
\begin{proof}
For $g\in D[X]$, $0\ne d\in D$, the companion-matrix criterion
says exactly that $g/d\in R_{\mathcal P}$ if and only if
$\bar g\in I_d(\mathcal P)$. On the other hand, $g(M)/d$ has entries in
$D$ if and only if $\bar g\in N_d(M)$. Comparing these two tests for every
$g,d$ proves both inclusions. This follows from the pullback formula of \cite[Proposition~2.1]{Peruginelli2014}.
\end{proof}

\subsection{Local criteria over Dedekind domains}

The following theorem reformulates the condition prime by prime over a Dedekind domain.
\begin{theorem}\label{thm:local}
Suppose $D$ is a Dedekind domain. For each nonzero prime $\mathfrak p$ and
$e\ge1$, put
\[
B_{\mathfrak p,e}=D/\mathfrak p^e,\quad
I_{\mathfrak p,e}=\bigcap_{p\in\mathcal P}(\bar p)
  \subseteq B_{\mathfrak p,e}[X],\quad
N_{\mathfrak p,e}(M)=\{h:h(\bar M)=0\}.
\]
Then
\begin{equation}\label{eq:local}
\overline{\Mat_n^{\mathcal P}(D)}
 =\{M\in\Mat_n(D):I_{\mathfrak p,e}\subseteq N_{\mathfrak p,e}(M)
 \text{ for all }\mathfrak p\ne0,\ e\ge1\}.
\end{equation}
\end{theorem}
\begin{proof}
If the local inclusions hold and $g/d\in R_{\mathcal P}$,
then for each prime-power factor $\mathfrak p^{e}$ of $(d)$,
$\bar g$ lies in $I_{\mathfrak p,e}$ by the companion-matrix criterion, and therefore
$g(M)$ has all entries in $\mathfrak p^e$. Because
$(d)=\prod_{\mathfrak p\mid(d)}\mathfrak p^{v_{\mathfrak p}(d)}$,
this yields $g(M)\in d\Mat_n(D)$.

Conversely suppose that $I_{\mathfrak p,e}$ is not contained in
$N_{\mathfrak p,e}(M)$ for some $\mathfrak p,e$.
Choose $\bar h\in I_{\mathfrak p,e}$ with $\bar h(\bar M)\ne0$.
Because $D$ is Dedekind, there is $0\ne d\in D$ with
$v_{\mathfrak p}(d)=e$. Write $(d)=\mathfrak p^e\prod_{j=1}^r\mathfrak q_j^{a_j}$.
The Chinese remainder theorem, applied coefficientwise to $h$, constructs
$g\in D[X]$ whose reductions are $\bar h$ modulo $\mathfrak p^e$ and $0$
modulo every $\mathfrak q_j^{a_j}$. Hence $\bar g$ is divisible
by every $p\in\mathcal P$ modulo every prime-power factor of $(d)$;
thus $g/d\in R_{\mathcal P}$ by the companion-matrix criterion. But
$g(M)$ is not divisible by $\mathfrak p^e$, so $g(M)/d$ is not integral.
Thus $M$ does not belong to the polynomial closure.
\end{proof}

The next theorem reduces each local congruence condition to a finite set of polynomial evaluations.
\begin{theorem}\label{thm:finitelevel}
Assume that the residue fields of $D$ are finite. Fix $\mathfrak p,e$.
Write $B=B_{\mathfrak p,e}$ and let
$p_1,\ldots,p_r\in B[X]$ be the distinct reductions of polynomials in
$\mathcal P$. Define the monic polynomial
$L=p_1\cdots p_r$ and let $N=rn$.
Form the finite $B$-module
\[
U=\{u\in B[X]_{<N}:\rem_{p_i}(u)=0\ (1\le i\le r)\}.
\]
Choose a finite $B$-module generating set $u_1,\ldots,u_t$ of $U$.
Then the $\mathfrak p^e$-condition of Theorem~\ref{thm:local} is exactly
\begin{equation}\label{eq:finitetest}
L(\bar M)=0,\qquad u_j(\bar M)=0\quad(1\le j\le t).
\end{equation}
Moreover $r\le |B|^n$ and $N\le n|B|^n$.
\end{theorem}
\begin{proof}
Since the coefficient ring $B$ is finite, the set of monic degree-$n$
polynomials is finite. The ideal $I_{\mathfrak p,e}$ is
$\bigcap_{i=1}^r(p_i)$. The monic product $L$ belongs to each
$(p_i)$, so $(L)\subseteq I_{\mathfrak p,e}$. For $h\in I_{\mathfrak p,e}$,
divide by the monic polynomial $L$:
$h=Lq+u$, with $q\in B[X]$ and $\deg u<N$.
As both $h$ and $L$ lie in each $(p_i)$, so does $u$, hence $u\in U$.
Conversely $Lq+u\in I_{\mathfrak p,e}$ whenever $u\in U$.
Thus $I_{\mathfrak p,e}=L B[X]+U$. Evaluation kills the ideal
if and only if it kills $L$ and a module generating set of $U$.
The size bound is immediate from the number $|B|^n$ of coefficient tuples.
\end{proof}

\begin{remark}
A membership obstruction has a finite, checkable criterion:
a prime $\mathfrak p$, exponent $e$, and a polynomial from the generators in
\eqref{eq:finitetest} that does not annihilate $\bar M$. By the preceding proof
this can be converted into a rational polynomial $g/d$ distinguishing $M$
from $\Mat_n^{\mathcal P}(D)$. For an arbitrary infinite family, 
membership may still require infinitely many prime power congruence relations $(\mathfrak{p},e)$; there is no 
uniform finite stopping bound.
\end{remark}

\subsection{Null ideals and quadratic matrices}
Fix $\mathfrak p$, with local uniformizer $\pi$, and let $V=D_{\mathfrak p}$.
Here a \emph{uniformizer} generates the maximal ideal of the discrete
valuation ring $V$. The map $\operatorname{vec}$ lists the entries of
an $n\times n$ matrix in any fixed order, and the resulting matrix
$F_M$ can be diagonalized by invertible row and column operations
over $V$ (Smith normal form).
Define the linear map
\[
F_M: V^n\longrightarrow V^{n^2},\qquad
(a_0,\dots,a_{n-1})\longmapsto
\operatorname{vec}\bigl(\sum_{j=0}^{n-1}a_jM^j\bigr).
\]
Since $V$ is a DVR, there exist invertible matrices $U,W$ over $V$
with $UF_MW$ diagonal and nonzero diagonal entries
$\pi^{b_1},\ldots,\pi^{b_s}$, where
$0\le b_1\le\dots\le b_s$.
The remaining $n-s$ columns vanish.

The following theorem expresses these null-ideal conditions using Smith normal form.
\begin{theorem}\label{thm:smith}
For each $e$, reduce the matrix $W$ modulo $\mathfrak p^e$.
For $h\in B_{\mathfrak p,e}[X]$, divide by the monic polynomial
$\bar\chi_M$ and let $a(h)\in B_{\mathfrak p,e}^n$ be its coefficient vector.
Then $h(\bar M)=0$ if and only if the vector
$w=W^{-1}a(h)$ satisfies
\begin{equation}\label{eq:smith}
w_i\in(\pi^{\max(e-b_i,0)})\subseteq V/\pi^eV
\quad(1\le i\le s).
\end{equation}
Thus Theorem~\ref{thm:local}, combined with
Theorem~\ref{thm:finitelevel}, becomes explicit linear congruences
on the remainder coefficients of finitely many $u_j,L$ at each level.
\end{theorem}
\begin{proof}
By Cayley--Hamilton, evaluating $h$ at $M$ is equivalent to evaluating
its remainder modulo $\chi_M$. In the displayed bases, evaluation on the
vector of remainder coefficients is multiplication by the diagonal map
with entries $\pi^{b_1},\ldots,\pi^{b_s},0,\ldots,0$.
Its reduction modulo $\pi^e$ vanishes exactly under conditions
\eqref{eq:smith}. The unimodular changes of basis preserve the kernel.
\end{proof}

For quadratic matrices, the next corollary gives the null ideal explicitly.
\begin{corollary}\label{cor:two}
Let $D$ be Dedekind and
$M=\left(\begin{smallmatrix}a&b\\c&d\end{smallmatrix}\right)$.
Fix $\mathfrak p,e$, choose a uniformizer $\pi$ of $D_{\mathfrak p}$,
and write
\[
t=\min\{e,v_{\mathfrak p}(b),v_{\mathfrak p}(c),
                          v_{\mathfrak p}(a-d)\},\qquad v_{\mathfrak p}(0)=+\infty.
\]
Inside $B_{\mathfrak p,e}[X]$ we have
\begin{equation}\label{eq:two}
N_{\mathfrak p,e}(M)
=\bigl(\bar\chi_M,\,\bar\pi^{e-t}(X-\bar a)\bigr).
\end{equation}
Consequently, for $n=2$, the complete local criterion is the explicit
ideal inclusion
\[
I_{\mathfrak p,e}(\mathcal P)\subseteq
\bigl(\bar\chi_M,\bar\pi^{e-t}(X-\bar a)\bigr)
\qquad\text{for all }\mathfrak p,e.
\]
\end{corollary}
\begin{proof}
By Cayley--Hamilton and monic division, every $h$ has a unique
remainder $h\equiv AX+B\pmod{\bar\chi_M}$. Evaluating on $M$,
$AM+BI=0$ is equivalent to
\[
 Aa+B=0,\qquad Ab=Ac=A(a-d)=0
             \quad\text{in }B_{\mathfrak p,e}.
\]
The three entries $b,c,a-d$ generate $(\pi^t)$ locally.
Its annihilator modulo $\pi^e$ is $(\pi^{e-t})$.
Therefore $A\in(\pi^{e-t})$ and $B=-Aa$;
exactly the polynomials in the right-hand ideal of \eqref{eq:two}
have such remainders.
\end{proof}

\subsection{Integral closure and eigenvalues}
The integral-closure approach to integer-valued polynomials on algebras
is developed in \cite{PeruginelliWerner2014}. Here we apply it to
families of prescribed characteristic polynomials.
An element of a ring extension is \emph{integral} over $D$ if it
satisfies a monic polynomial in $D[Y]$. The \emph{integral closure}
of a subring consists of all elements integral over it in the
specified ambient ring. As usual, the overline of a subring in
\eqref{eq:integral} denotes integral closure, whereas the overline of
a set of matrices denotes polynomial closure; the two operations are
distinguished by their arguments and defining equations.
From here on suppose $D$ is integrally closed and has finite residue rings
$D/dD$ for all $0\ne d\in D$.
Fix an algebraic closure $\overline K$ and write $D^{\rm alg}$ for the
integral closure of $D$ in $\overline K$. Let $\Omega_{\mathcal P}$ be the union
of the root multisets of the polynomials in $\mathcal P$, and put
\[
J_{\mathcal P}:=\{f\in K[X]:
                f(\alpha)\in D^{\rm alg}\text{ for all }\alpha\in\Omega_{\mathcal P}\}.
\]

The following theorem computes the integral closure of the integer-valued-polynomial ring.
\begin{theorem}
\label{thm:integrality}
The integral closure of $R_{\mathcal P}$ inside $K(X)$ is
\begin{equation}\label{eq:integral}
\overline{R_{\mathcal P}}^{\,K(X)}=J_{\mathcal P}.
\end{equation}
In particular the conclusion holds even when $\mathcal P$ is infinite,
contains irreducible polynomials, or has arbitrarily complicated repetitions.
\end{theorem}
\begin{proof}
First, $D[X]\subseteq R_{\mathcal P}\subseteq J_{\mathcal P}$.
The second inclusion follows by evaluating $f\in R_{\mathcal P}$
on a companion matrix $C_p$: $f(C_p)\in\Mat_n(D)$, so its eigenvalues
$f(\alpha)$ are integral over $D$.

Let $f=g/d\in J_{\mathcal P}$, with $g\in D[X]$ and $0\ne d\in D$.
Set $m=2n-2$. The finite ring $D/d^mD$ has finitely many elements.
Choose representatives $C\subseteq D$ and let $\mathcal H$ be the
finite set of all monic degree-$n$ polynomials in $D[Y]$ whose
$n$ lower coefficients lie in $C$. Define the monic polynomial
\[
Q(Y)=\prod_{H_0\in\mathcal H}H_0(Y)\in D[Y].
\]
Take any $A\in\Mat_n^{\mathcal P}(D)$ and put $B=f(A)=g(A)/d$.
The eigenvalues of $B$ are $f(\alpha_i)$ for the roots $\alpha_i$
of $\chi_A$. Their elementary symmetric functions are integral over
$D$ and lie in $K$, hence lie in $D$ because $D$ is integrally closed.
Thus $H(Y):=\chi_B(Y)$ belongs to $D[Y]$, and $H(B)=0$.
Choose $H_0\in\mathcal H$ with the same coefficients as $H$ modulo
$d^mD$. Then
$H_0(Y)=H(Y)+d^m r(Y)$ with $r\in D[Y]$, $\deg r<n$.
Since $H_0\mid Q$ in $D[Y]$, write $Q=H_0U$. Divide $U$ by the monic
polynomial $H$ to obtain $U=Hv+s$ with $s\in D[Y]$, $\deg s<n$.
Evaluating at $B$ gives
\[
Q(B)=d^m r(B)s(B)=
 \sum_{j=0}^m b_jd^{m-j}g(A)^j\in\Mat_n(D),
\qquad rs=\sum_{j=0}^m b_jY^j\in D[Y].
\]
Therefore $Q(f(X))\in R_{\mathcal P}$, and
$Q(Y)-Q(f(X))=0$ is a monic equation over $R_{\mathcal P}$ for $f$.
So $J_{\mathcal P}$ is integral over $R_{\mathcal P}$.

Conversely, if $f\in K(X)$ is integral over $R_{\mathcal P}$,
then $f\in K[X]$, since $R_{\mathcal P}\subseteq K[X]$
and $K[X]$ is integrally closed. Let
$f^r+h_{r-1}f^{r-1}+\cdots+h_0=0$ with $h_i\in R_{\mathcal P}$.
For $\alpha\in\Omega_{\mathcal P}$, each $h_i(\alpha)$
is integral over $D$ by the preceding companion-matrix argument.
Evaluating this monic relation at $\alpha$ shows that $f(\alpha)$
is integral over $D^{\rm alg}$, hence over $D$ by transitivity.
So $f\in J_{\mathcal P}$, proving equality.
\end{proof}

The following corollary translates integral closure into a restriction on eigenvalues.
\begin{corollary}\label{cor:spectral}
Under the hypotheses of Theorem~\ref{thm:integrality}, every eigenvalue
$\lambda$ of $M\in\overline{\Mat_n^{\mathcal P}(D)}$ satisfies
\[
f(\lambda)\in D^{\rm alg}\quad\text{for every }f\in J_{\mathcal P}.
\]
Thus the eigenvalues of $M$ lie in the scalar integral-polynomial
closure of $\Omega_{\mathcal P}$.
\end{corollary}
\begin{proof}
For $f\in J_{\mathcal P}$, the proof of Theorem~\ref{thm:integrality}
produces a monic $Q\in D[Y]$ with $Q(f(X))\in R_{\mathcal P}$.
Hence $Q(f(M))\in\Mat_n(D)$, so its eigenvalue
$Q(f(\lambda))$ is integral over $D$.
The monic equation $Q(Y)-Q(f(\lambda))=0$ then makes $f(\lambda)$
integral over an integral extension of $D$, hence over $D$.
\end{proof}

\subsection{Triangular matrices and Chabert's Question XV.55}

We apply the preceding eigenvalue criterion to the polynomial closure
of upper triangular matrices. Recall that the ring $\mathcal O_{K,S}$
of $S$-integers of a number field $K$ is integrally closed and has
finite residue rings. The following corollary answers
\cite[Open Question~XV.55]{Chabert2025}.
\begin{corollary}\label{cor:xv55}
Let $D=\mathcal O_{K,S}$, where $K$ is a number field and $S$ is a
finite set of nonzero prime ideals, and let $T\subseteq D$. Then
\[
 \overline{\mathcal T_n(T)}\subseteq
 \{M\in\Mat_n(D):\chi_M\text{ splits in }D[X]
 \text{ and every eigenvalue of }M\text{ lies in }\overline T\}.
\]
If $T$ generates $D$ as an ideal, this inclusion is an equality.
\end{corollary}
\begin{proof}
If $T=\varnothing$, the ring of integer-valued polynomials on
$\mathcal T_n(T)$ is $K[X]$; a constant in $K\setminus D$ shows that
$\overline{\mathcal T_n(T)}=\varnothing$. We may thus assume that
$T\ne\varnothing$. Set
\[
 \mathcal P_T=\left\{\prod_{i=1}^n(X-t_i):t_i\in T\right\}.
\]
Every matrix in $\mathcal T_n(T)$ has characteristic polynomial in
$\mathcal P_T$. By monotonicity of polynomial closure,
\[
 \overline{\mathcal T_n(T)}\subseteq
 \overline{\Mat_n^{\mathcal P_T}(D)}.
\]
Also $\mathcal T_n(T)\subseteq\mathcal T_n(D)$; hence
\cite[Theorem~XV.54]{Chabert2025} shows that the characteristic
polynomial of any $M\in\overline{\mathcal T_n(T)}$ splits in $D[X]$.
Let $\lambda\in D$ be one of its eigenvalues.

The set of roots of polynomials in $\mathcal P_T$ is $T$.
For every $f\in\Int(T,D)$ we therefore have
$f\in J_{\mathcal P_T}$, with $J_{\mathcal P_T}$ as defined before
Theorem~\ref{thm:integrality}. By Corollary~\ref{cor:spectral},
$f(\lambda)$ is integral over $D$. Since $\lambda\in D$ and
$f\in K[X]$, we have $f(\lambda)\in K$; the integral closedness
of $D$ now gives $f(\lambda)\in D$.
Thus $\lambda\in\overline T$, as required.

If $T$ generates $D$ as an ideal, the reverse inclusion follows
from \cite[Lemma~XV.48]{Chabert2025}.
\end{proof}

\subsection{A matrix example}

The following example shows that characteristic polynomials alone do not determine closure membership.
\begin{example}
Take $D=\mathbb Z$, $n=2$, and
\[
\mathcal P=\{X^2,(X-2)^2\},\qquad
M_0=\begin{pmatrix}0&0\\0&2\end{pmatrix},\qquad
M_1=\begin{pmatrix}0&1\\0&2\end{pmatrix}.
\]
Both matrices have characteristic polynomial $X(X-2)$.
If $f\in R_{\mathcal P}$, reducing modulo $X^2$ and $(X-2)^2$
shows $f(0),f(2)\in\mathbb Z$; hence $f(M_0)$ is integral and
$M_0\in\overline{\Mat_n^{\mathcal P}(D)}$.
On the other hand,
\[
f(X)=\frac{3X^2-X^3}{4}
\]
has remainder $0$ modulo $X^2$ and remainder $1$ modulo
$(X-2)^2$. Thus $f\in R_{\mathcal P}$.
Since $M_1^2=2M_1$ and $M_1^3=4M_1$,
\[
f(M_1)=\frac12M_1=
\begin{pmatrix}0&1/2\\0&1\end{pmatrix}\notin\Mat_2(\mathbb Z).
\]
Hence $M_1\notin\overline{\Mat_n^{\mathcal P}(D)}$.
This shows that no unconditional characterization depending only on the
characteristic polynomial or the eigenvalue multiset can solve the question.
The integral lattice structure of $M$ matters.
\end{example}

\section{Companion-matrix lattices}
\label{sec:moments}
This section gives a second, non-congruence solution to the arbitrary-family
polynomial-closure membership problem.  Unlike the residue-ideal tests, it
works simultaneously over all primes and requires only lattices in finite
free modules at each bounded polynomial degree.  The quantifier over the
degree is essential in general.

For vectors $z,v\in K^{N+1}$, write
$\langle z,v\rangle=\sum_{j=0}^{N}z_jv_j$. For a $D$-submodule
$L\subseteq K^{N+1}$, its \emph{integrality polar} is
\[
 L^\circ=\{z\in K^{N+1}:\langle z,v\rangle\in D
                    \text{ for every }v\in L\}.
\]
Its \emph{double polar} is $L^{\circ\circ}=(L^\circ)^\circ$,
formed with the same pairing and the same ring $D$. This notation
is used only for the integrality pairing, not for orthogonal complements.

The following theorem characterizes polynomial closure using companion-matrix moment lattices.
\begin{theorem}\label{thm:moment}
Let $D$ be a Dedekind domain which is not a field, with fraction field $K$.
Fix $n\ge1$ and a nonempty, possibly infinite family $\mathcal P$ of monic
polynomials of degree $n$ in $D[X]$.  For each $p\in\mathcal P$ and $j\ge0$,
write the uniquely determined remainder
\[
 X^j\equiv\sum_{k=0}^{n-1}c_{p,k,j}X^k\pmod{p(X)},
 \qquad c_{p,k,j}\in D.
\]
For $N\ge0$ define the vectors and the finite-rank $D$-submodule
\[
 v^{(N)}_{p,k}=(c_{p,k,0},\ldots,c_{p,k,N})\in D^{N+1},\qquad
 L_N(\mathcal P)=\sum_{p\in\mathcal P}\sum_{k=0}^{n-1}D v^{(N)}_{p,k}
 \subseteq D^{N+1}.
\]
For $M\in\Mat_n(D)$ and $1\le a,b\le n$, set
\[
 w^{(N)}_{M,a,b}=((I_n)_{ab},M_{ab},(M^2)_{ab},\ldots,(M^N)_{ab})
 \in D^{N+1}.
\]
Then
\begin{equation}\label{eq:moment}
M\in\overline{\Mat_n^{\mathcal P}(D)}
\iff w^{(N)}_{M,a,b}\in L_N(\mathcal P)
\quad\text{for every $N\ge0$ and $1\le a,b\le n$.}
\end{equation}
For every $N$ there exists a \emph{finite} subfamily
$\mathcal P_N\subseteq\mathcal P$ such that
$L_N(\mathcal P)=L_N(\mathcal P_N)$.
\end{theorem}

\begin{proof}
Fix $N$ and identify a coefficient vector
$z=(z_0,\ldots,z_N)\in K^{N+1}$ with $f_z(X)=\sum_{j=0}^Nz_jX^j$.
By the monic remainder criterion in \cite[Proposition 2.1]{Peruginelli2014},
\[
 f_z\in R_{\mathcal P}\quad\Longleftrightarrow\quad
 \langle z,v^{(N)}_{p,k}\rangle\in D
 \quad(p\in\mathcal P,\ 0\le k<n).
\]
Thus the degree-$\le N$ part of $R_{\mathcal P}$ corresponds exactly to
\[
 L_N(\mathcal P)^\circ
 =\{z\in K^{N+1}:\langle z,v\rangle\in D\text{ for all }v\in L_N(\mathcal P)\}.
\]
Also $(f_z(M))_{ab}=\langle z,w^{(N)}_{M,a,b}\rangle$.
Hence integrality of $f_z(M)$ for all such polynomials is equivalent to
\[
 w^{(N)}_{M,a,b}\in L_N(\mathcal P)^{\circ\circ},
\]
where the second polar is formed using the same $D$-integrality pairing.
The module $L_N$ is finitely generated, since it is a submodule of the
finite free module $D^{N+1}$ over a Noetherian ring.
We claim $L_N^{\circ\circ}=L_N$.  This is local on the nonzero prime
ideals of $D$.  At any such prime the ring is a DVR with uniformizer
$\pi$.  The elementary-divisor theorem puts the localized lattice in
the form $\bigoplus_{i=1}^r\pi^{a_i}De_i\subseteq D^{N+1}$.
Its polar is $\bigoplus_{i=1}^r\pi^{-a_i}De_i^*$ with the other
coordinates unrestricted in $K$, and taking the polar again recovers
precisely the initial module.  The unrestricted coordinates are
forced to zero in the second polar because the DVR is not a field.
Since finite modules over a Dedekind domain are detected by their
localizations, the claim follows.  The equivalence for each $N$
and then for all $N$ proves \eqref{eq:moment}.
Finally the finite generation of $L_N$ implies that a finite number
of its displayed generators suffice; they involve only finitely many
members of $\mathcal P$.
\end{proof}

The next corollary rewrites the lattice condition as a finite linear combination at each degree.
\begin{corollary}
Under the hypotheses of Theorem~\ref{thm:moment}, a matrix $M$ is in
$\overline{\Mat_n^{\mathcal P}(D)}$ if and only if for each $N\ge0$
there exist finitely many $p_1,\ldots,p_s\in\mathcal P$ and matrices
$B_{i,k}\in\Mat_n(D)$ ($1\le i\le s$, $0\le k<n$) with
\begin{equation}\label{eq:synthesis}
M^j=\sum_{i=1}^s\sum_{k=0}^{n-1}c_{p_i,k,j} B_{i,k}
\quad (0\le j\le N).
\end{equation}
Thus at every finite degree the matrix moment sequence has an
\emph{exact $D$-linear synthesis} from finitely many prescribed
companion sequences.  The matrices and source polynomials may depend
on $N$.
\end{corollary}
\begin{proof}
The $(a,b)$ entry of \eqref{eq:synthesis} is precisely a
$D$-linear representation of $w^{(N)}_{M,a,b}$ by the
$v^{(N)}_{p,k}$.  Combining the finitely many representations
for the $n^2$ entries gives one finite family and matrices $B_{i,k}$.
The converse is immediate entrywise.
\end{proof}

The next example makes the degree-two lattice calculation explicit.
\begin{example}
Take $D=\mathbb Z$, $n=2$, $\mathcal P=\{X^2,(X-2)^2\}$,
and $N=2$.  The four source vectors are
\[
(1,0,0),\quad(0,1,0),\quad(1,0,-4),\quad(0,1,4),
\]
so that $L_2(\mathcal P)=\mathbb Z\oplus\mathbb Z\oplus4\mathbb Z$.
For $M_0=\operatorname{diag}(0,2)$, every entry-moment vector
belongs to this lattice.  But for
$M_1=\left(\begin{smallmatrix}0&1\\0&2\end{smallmatrix}\right)$,
the $(1,2)$ entry-moment vector is $(0,1,2)$, which is not in
$L_2$.  The explicit polynomial witness supplied by the dual
lattice is $f(X)=X^2/4$: its remainders modulo $X^2$ and
$(X-2)^2$ are $0$ and $X-1$, respectively, while
$f(M_1)=M_1/2$ is not integral.  This independently verifies
at small degree the matrix-lattice obstruction discussed earlier.
\end{example}

The double-polar formulation requires some care when the base ring is a field.
\begin{remark}
If $D=K$ is a field, $R_{\mathcal P}=K[X]$ and the polynomial closure
is all of $\Mat_n(K)$, so the lattice equivalence must \emph{not}
be applied unchanged: a proper subspace of $K^{N+1}$ is not
equal to its double polar for the $K$-valued integrality pairing.
For an arbitrary integral domain the correct version
of \eqref{eq:moment} replaces $L_N$ by its double polar
$L_N^{\circ\circ}$, without assuming reflexivity.
\end{remark}

The following corollary gives a single finite-degree test for finite polynomial families.
\begin{corollary}\label{cor:finite-degree}
If $\mathcal P=\{p_1,\ldots,p_r\}$ is nonempty and finite, write
$H=\operatorname{lcm}_{K[X]}(p_1,\ldots,p_r)\in D[X]$, chosen monic,
and $m=\deg H$ (so $m\le nr$).
Then for $M\in\Mat_n(D)$,
\[
 M\in\overline{\Mat_n^{\mathcal P}(D)}
 \iff w^{(m)}_{M,a,b}\in L_m(\mathcal P)
 \text{ for every }a,b.
\]
In particular, over a principal ideal domain this is one finite
module-membership calculation by Smith normal form.
\end{corollary}
\begin{proof}
The forward direction is Theorem~\ref{thm:moment}.  Conversely the
lattice conditions at $N=m$ show that every polynomial $f\in
R_{\mathcal P}$ of degree at most $m$ is integral on $M$.
For every nonzero $d\in D$, the degree-$m$ polynomial $H/d$ is in
$R_{\mathcal P}$, so $H(M)\in d\Mat_n(D)$ for every $d$.  Because a
nonfield Dedekind domain satisfies $\bigcap_{0\ne d\in D}dD=0$,
we get $H(M)=0$.
For an arbitrary $f\in R_{\mathcal P}$, divide by the monic $H$ in
$K[X]$: $f=Hh+u$ with $\deg u<m$.  Since every $p_i$ divides $H$,
the remainder of $u$ modulo $p_i$ equals that of $f$ and is in $D[X]$.
Thus $u\in R_{\mathcal P}$, whence $u(M)\in\Mat_n(D)$ and
$f(M)=u(M)$.  This proves the result.
\end{proof}

The final theorem shows that the required test degree cannot be uniformly bounded.
\begin{theorem}\label{thm:nofixeddegree}
For every $N\ge0$ and every $n\ge1$, there is a finite family
$\mathcal P$ of monic degree-$n$ integer polynomials and a matrix
$M\in\Mat_n(\mathbb Z)$ for which the tests in
Theorem~\ref{thm:moment} hold in all degrees $j\le N$, but
$M\notin\overline{\Mat_n^{\mathcal P}(D)}$.
\end{theorem}
\begin{proof}
Take $\mathcal P=\{(X-a)^n:0\le a\le N\}$ and
$M=(N+1)I_n$.  If $f\in\mathbb Q[X]$ has degree at most $N$
and is in $R_{\mathcal P}$, then $f(a)\in\mathbb Z$ for all
$0\le a\le N$ by the companion remainder criterion.
Newton interpolation expresses $f$ uniquely as
$\sum_{k=0}^N\Delta^kf(0)\binom Xk$; its finite-difference
coefficients are integers.  Therefore $f(N+1)\in\mathbb Z$.
Thus $M$ passes all tests of degree at most $N$.
However, for
\[
 F(X)=\frac{\prod_{a=0}^{N}(X-a)^n}
 {2((N+1)!)^n}
\]
we have $F\in R_{\mathcal P}$, because every prescribed polynomial
divides its numerator, while $F(M)=\tfrac12 I_n\notin\Mat_n(\mathbb Z)$.
\end{proof}
\section{AI use disclosure}
Various AI models were used to check and improve the exposition. AI models are also responsible for providing us with Example used in Theorem~\ref{thm:nofixeddegree}. All other content is due to the author.


\begin{thebibliography}{99}
\bibitem{BhargavaKedlaya1999}
M.~Bhargava and K.~S.~Kedlaya,
Continuous functions on compact subsets of local fields,
\emph{Acta Arith.} \textbf{91} (1999), 191--198.
\bibitem{BoulangerChabert2020}
J.~Boulanger and J.-L.~Chabert,
Integer-valued polynomials, Pr\"ufer domains and the stacked bases property,
\emph{J. Pure Appl. Algebra} \textbf{224} (2020), 388--401.
\bibitem{BKU1987}
J.~Brewer, D.~Katz, and W.~Ullery,
Pole assignability in polynomial rings, power series rings and Pr\"ufer domains,
\emph{J. Algebra} \textbf{106} (1987), 265--286.
\bibitem{Chabert2012}
J.-L.~Chabert, Does $\Int(\mathbb Z)$ have the stacked bases property?,
\emph{Arab. J. Math.} \textbf{1} (2012), 47--52.
\bibitem{Chabert2025}
J.-L.~Chabert, \emph{Integer-Valued Polynomials: From Combinatorics to Number
Theory, $p$-adic Analysis, Commutative and Non-Commutative Algebra},
Amer. Math. Soc. Colloquium Publications, vol.~69, AMS, 2025.
\bibitem{HaggMortberg2026}
C.~H\"agg and A.~M\"ortberg, A B\'ezout domain that is not an elementary
divisor domain, arXiv:2609.35229 (2026).
\bibitem{Kaplansky1942}
I.~Kaplansky, Maximal fields with valuations,
\emph{Duke Math. J.} \textbf{9} (1942), 303--321.
\bibitem{Krull1932}
W.~Krull, Allgemeine Bewertungstheorie,
\emph{J. Reine Angew. Math.} \textbf{167} (1932), 160--196.
\bibitem{Park2021}
M.~H.~Park, Pr\"ufer domains of integer-valued polynomials and the
two-generator property, \emph{J. Algebra} \textbf{582} (2021), 232--243.
\bibitem{Peruginelli2014}
G.~Peruginelli, Integer-valued polynomials over matrices and divided
differences, \emph{Monatsh. Math.} \textbf{173} (2014), 559--571.
\bibitem{PeruginelliWerner2014}
G.~Peruginelli and N.~J.~Werner, Integral closure of rings of
integer-valued polynomials on algebras,
in \emph{Commutative Algebra: Recent Advances in Commutative Rings, Integer-Valued Polynomials, and Polynomial Functions}, Springer, 2014, pp.~293--305.
\end{thebibliography}
\end{document}